  \RequirePackage{lineno} 
\documentclass[12pt]{article}
\usepackage{amssymb}

\usepackage{amssymb,amsfonts}

\usepackage{pdflscape}
\usepackage{graphicx}
\usepackage{mathtools}
\usepackage{amsthm}
\usepackage{setspace}
\usepackage{diagbox}
\usepackage{adjustbox}
\usepackage{array}
\usepackage{multirow}
\usepackage{booktabs}
\usepackage[all,arc]{xy}
\usepackage{enumerate}
\usepackage{mathrsfs}
\usepackage{setspace}
\newtheorem{definition}{Definition}

\theoremstyle{definition}

\theoremstyle{remark}

\makeatletter
\let\c@equation\c@thm
\makeatother
\numberwithin{equation}{section}

\title{Formation of  coalition structures as a non-cooperative game 2: applications }

\author{Dmitry Levando, \thanks{ 
Special thanks to Fuad Aleskerov,   Shlomo Weber and Lev Gelman.
 \newline E-mail for correspondence: dlevando ({at}) hse.ru.}}
\date{} 

\begin{document}

\maketitle

\begin{abstract}
The paper uses a non-cooperative simultaneous   game  for  coalition structure formation  (Levando, 2016) to demonstrate some   applications of the   introduced game: a cooperation, a Bayesian game within a coalition with   intra-coalition externalities, a stochastic game, where states are coalition structures, self-enforcement properties  of non-cooperative equilibrium and construction of  a non-cooperative stability criterion. 
 \end{abstract}



 



\section{Subject of the paper}
The previous paper, Levando (2016), introduced a non-cooperative game to study coalition structures formation as a non-cooperative game. The suggested game consisted of a set of players $N$\footnote{Size of the set $N$ will be also denoted as $N$, if there is no ambiguity in notation.}, a coalition structure construction mechanism ($\{ K, \mathcal{P}(K), \mathcal{R} (K)\}$), and   individual properties of players, $(S_i (K), \mathcal{U}(K))_{i \in N}$, such that $$\Gamma(K)=\Big\langle  K, \{ K, \mathcal{P}(K), \mathcal{R} (K)\}, (S_i (K), \mathcal{U}_i(K))_{i \in N}\Big\rangle, $$
where $K$  is a maximum coalition size or a  maximum number of deviators, $\mathcal{P}(K)$ is a  family of coalition structures such that any coalition size (a number of deviators ) does not exceed $K$,  $\mathcal{R} (K)$ is a family of  specific coalition structure formation rules,
$S_i (K)$ is a corresponding set of strategies of $i \in N$ induced by $K$, $\mathcal{U}_i(K)$ is family of specific coalition structure  payoffs for $i$ in every  partition $P \in \mathcal{P}(K)$.
Games for all relevant $K$  make a nested   family $$\Gamma  = \{\Gamma(K=1), \ldots, \Gamma(K), \ldots \Gamma(K=N) \} \colon \Gamma(K=1) \subset  \Gamma(K) \subset \Gamma(K=N).$$

The family of games is characterized by a list of  equilibrium strategy profiles, $\Big( \sigma^*(1), \ldots, \sigma^*(K), \ldots, \sigma^*(K=N)\Big)$, where  $ \sigma^*(K)=( \sigma^*_i(K))_{i \in N}$. and by a list of  equilibrium  partitions,
 $\Big(\{ P^*\}(1), \ldots, \{ P^*\}(K),\ldots, \{ P^*\}(K=N) \Big) $, $\{ P^*\}(K)  \subset \mathcal{P}(K)$. Every game $\Gamma (K)$  from a family $\Gamma$, which  has an equilibrium may be in mixed strategies.

This paper demonstrates   applications of  the constructed  game  for: self-interest cooperation, a Bayesian stochastic game, studying  self-enforcing properties of an equilibrium and construction of a  non-cooperative criterion of coalition structure stability.

\section{Cooperation}
To explain formally cooperation we will retell  with necessary changes the example from the corporate dinner game in Levando (2016). The result of this example will be used to explain    formal definition  of  what   a cooperation in allocations of self-interest     players can be. In this section  the term "cooperation" means only a cooperation of players in allocations over coalitions. Cooperation in payoffs 
 is not addressed here.

 \subsection{Cooperation in a corporate dinner game}
 
 Consider a game of 4 players: A  is  a President; B is a senior vice-president; $C_1,C_2$ are  two other vice-presidents. 
Coalition  is  a group of players   at \textit{one} table. Every player may sit only at one table.  Coalition structure is an  allocation of  all players  over no more than four tables. Empty tables are not   taken into account.

 Individual set of strategies is  a set of all coalition structures  for the   players. Any  player can choose a coalition structure. 
  A set of  strategies in the game is a direct (Cartesian) product of four individual strategy sets.  A choice of all players is a point in the set of strategies of the game. 
  
Preferences are  such that  everyone (besides A) would like   to have a dinner with A, but A only with B.  Everyone wants players outside his table to eat individually,  due to possible dissipation of rumors or information exchange. No one can   enforce  others to form or not to form coalitions.

  In  every partition any   coalition (i.e. a table with players  who eat at the table) is formed only if everybody at the table agrees to have dinner together, otherwise the player  eats alone. The same coalition  belongs to different coalitions, but with different allocations of other players over  other coalitions. 
 
  The game is simultaneous and one shot. Realization  of the  final partition ( a coalition structure) depends on choices of coalition structures of  all players. 
  Example. Let player A choose $\{\{A,B \}, \{ C_1\},  \{ C_2\} \}$;  player  B choose $\{\{A,B \}, \{ C_1\},  \{ C_2\} \}$; player  $C_1$ choose $\{\{A,C_1 \}, \{ B\},  \{ C_2\} \}$,      and player  $C_2$ choose $\{\{A,C_2 \}, \{B\},  \{ C_2\} \}$. Then the final partition is  $\{\{A,B \}, \{ C_1\},  \{ C_2\} \}$.  It is clear that  strong Nash equilibrium (Aumann, 1960), which is based on  a deviation of a coalition of any size, does not  discriminate between coalition structures mentioned above.
  
Players have preferences over coalition structures. Payoff profile of the game should be defined for every final coalition structure. \
Table \ref{dinner} presents coalition structures only  with the best individual payoffs.\footnote{All other coalition structures have significantly  lower payoffs. } Thus only some partitions from the big set of all strategies deserve attention. The first column is  a number of a strategy. The second column is an allocation of players over a coalition structure, i.e. a strategy, and also  a list of   best final coalition structures. The third column is  a  payoff profile for all players of every listed coalition structure. The fourth column is a list of coalition  values in  the listed coalition structure, if to calculate values from  the cooperative game theory. 

Payoffs are defined for a game of four players in such a way, that an increase in the size of a maximum coalition from $K=2$ to $K=3,4$  only  decreases individual payoffs. Thus from one side there is "an optimal" coalition size, from another an increase in a maximum coalition size does not change the equilibrium. The last property  will be addressed formally in one of the following sections of the current  paper.

  \begin{table}[htp]
\caption{Strategies and payoffs in the corporate dinner game} 
\begin{center}
\begin{tabular}{|c|c|c|c|}
\hline 
{\small{num}} & {\small{Best final  partitions}} & 
\begin{tabular}{c}  
{\small{ Non-cooperative payoff profile }} 
\\
{\small{ $(U_A,U_B,U_{C_1},U_{C_2})$ }} 
\end{tabular}
   & \begin{tabular}{c}  
{\small{   Values of coalitions as in  }} 
\\
{\small{  cooperative  game theory}} 
\end{tabular} \\
\hline
1&  $\{ \{ A,B\}, \{ C_1\}, \{ C_2\}  \}$ & (10,10,3,3)  & $20_{AB}$, $3_{C_1}$, $3_{C_2}$\\
 $2^*$& $\{ \{ A,B\}, \{ C_1, C_2\}  \}$ & (8,8,5,5)  & $16_{AB}$, $10_{C_1,C_2}$\\
 3&$\{ \{ A,C_1\}, \{ B, C_2\}  \}$ & (3,5,10,5)  & $13_{AC_1}$, $10_{BC_2}$ \\
  4&$\{ \{ A,C_1\}, \{ B\} , \{C_2\}  \}$ & (3,3,10,3) & $13_{AC_1}$, $3_B$, $3_{C_2}$ \\
 5&$\{ \{ A,C_2\}, \{ B, C_1\}  \}$ & (3,5,5,10)  &  $13_{AC_2}$, $10_{BC_1}$   \\
 6&  $\{ \{ A,C_2\}, \{ B\} , \{C_1\}  \}$ & (3,3,3,10) & $13_{AC_1}$, $3_B$, $3_{C_1}$\\
\hline 7& \small{all other partitions} & (1,1,1,1)  & $\le 4 $\\
\hline
\end{tabular}
\end{center}
\label{dinner}
\end{table}%
The game runs as follows. Players simultaneously and independently announce choices of   desirable coalition structures, then the final coalition structure  is formed according to the rule above, and payoffs are assigned. The rule is not formalized  for a shorter exposition.

Players A and B would always like  to be    together, then with $C_1$ or $C_2$. Being rational  $A$ and $B$  would choose   coalition structure with the highest payoff, i.e. the strategy 1. The first best choices of $C_1$ and $C_2$  is    to have a dinner   with  $A$. However, A will never choose to be with either of them. Hence,  non-realization of this option makes a potential loss for both $C_1$ and $C_2$. 

Both players  $C_1$ and $C_2$ are aware of this. 
Another common knowledge,   that they both know about each other, is that  if they make a coalition of two, each will be better off. 
This behavior meets sociological understanding of cooperation - to unite in order to prevent a  loss, to which each is individually exposed to. 

 A cooperation between $C_1$ and $C_2$  does not demolish  the coalition $\{A,B \}$, but  only  reduces payoffs  for it{'}s members.  At the same time   players A and B cannot  prevent cooperation between $C_1$ and $C_2$ (or to insure against). 

On the other hand, if players A and B choose strategy 2  they will  obtain coalition $\{A,B \}$ in any case.  
 In terms of mixed strategies this means that an equilibrium mixed strategy for $A$ and for $B$ is the whole probability space over two points, two coalition structures.
 
From the forth column we can see that the corresponding cooperative game has an empty core. Strong Nash equilibrium cannot be applied here also. Coalition value approach also supports the idea of cooperation, but without an explicit allocation of payoffs. This means it ignores individual rationality and incentive compatibility to participate in a coalition. For example, how players $C_1$ and $C_2$ should divide the value of a joint coalition $\{ C_1,C_2\}$ equal to  10?

The constructed game has a unique equilibrium,  which in terms of individual payoffs is characterized as the second-best efficient for every player. There are   two coalitions in the equilibrium coalition structure.


\subsection{Formal definition of cooperation}

This section formalizes an idea of cooperation presented in the example above. We demonstrate only  one way for defining cooperation, when it is intentional and  complete.



\begin{definition}[complete cooperation in a coalition] 
In a game $\Gamma(K)$ with an equilibrium $\sigma^*(K)=(\sigma^*_i)_{i \in N}$ a set of players $g$, \textbf{ completely cooperate  in   a coalition   $g$  }  \textit{ex ante} if for every player $i \in g$ there is 
\begin{description}
\item[\textit{ex ante}]: for every $i$ in $g$, a desirable coalition $g$ always belongs to a chosen coalition structure, i.e such if $s_i$ is chosen by $i$, then $s_i \in S_i(P_i)$,  $g \in P_i$, where $P_i$ is a coalition structure chosen by $i$, however coalition structures chosen by different players may be different.   
\item[\textit{ex post 1}]: every realized coalition structure contains $g$, i.e.  $g \in  \forall P^*$, where $P^*$ is a formed equilibrium partition of $\Gamma(K)$,
\item[\textit{ex post 2}]   $EU^{\Gamma(K)}_i \Big(\sigma^*(K)\Big) \ge EU^{\Gamma(K)}_i \Big(\sigma^*_i(K), \sigma^*_{-i}(K)\Big)$, for  $\forall \sigma_i(K) \neq \sigma^*_i(K)$,
\end{description}
\end{definition}

First of all,   cooperation is defined for a game $\Gamma(K)$.  If a game changes due to a change in the number of deviators $K$, the cooperation may evaporate. 

Every player chooses a partition with the desirable coalition $g$ with positive probability. 
 But individually chosen coalition structures by all players may be  different.  
 
 After the game is over the coalition $g$ always belongs to every  final equilibrium coalition structure, disregard allocation of players over other coalitions.  A final partition may  differ from a chosen one, but in any case it will contain the desirable coalition.
 
 The defined  
 cooperation assumes agents are acting from  self-interest motivations. The  lunch game example further expands  the case, where there is  imperfect cooperation.

The  dinner game example  has the  complete   cooperation for players $C_1$  and $C_2$. 
The definition does not exclude inter and intra-coalition interaction.  
If we relax some conditions of the definition we will obtain weaker conditions for cooperation.  Cooperation in repeated games is of special interest and will be addressed in one  of the next papers.


\section{Bayesian games}
In this section we demonstrate how     intra-coalition  externalities of the grand coalition may happen from equilibrium mixed strategies. In order to show  that, a standard Battle of Sexes game is modified.

Let there be two players, Ann and Bob. Each  has  two options:  to be together with another or to be alone. In every option each can choose   where to go: to Box or to Opera. Hence every player has four strategies. A set of strategies  of the game leads to 16 outcomes. Every outcome consists of payoff profile and a partition (or a coalition structure). Both players have preferences over coalition structures: they prefer to be together, then be separated. 

The rules of coalition structure formation mechanism are:
\begin{enumerate}
\item If they both choose to be together, i.e. both choose  the coalition structure $P_{joint}=\{Ann, Bob \}$ then:
\begin{enumerate}
\item if both choose the same  action for $P_{joint}$ (i.e. both choose Box or both choose Opera), then they go to where they both choose  to go,
\item otherwise they  do not go anywhere, but enjoy just being together;
\end{enumerate}
\item if at least one of them chooses to stay alone, i.e. chooses a partition   $P_{separ}=\{\{Ann\},\{ Bob \}\}$, then each goes alone to where she/he chooses, may be to different Opera or different Box performances.
\end{enumerate}

  {\tiny{\begin{table}[htp]
\caption{Payoff for  the  Bach-or-Stravinsky game.  {B is for Box, O is for Opera. If the players choose to be together, and it is realized due to the rule of coalition structure formation, then   each obtains   an additional fixed payoff $\epsilon >0$.}
} 
\begin{center}
\begin{tabular}{|c|c|c||c|c|} \hline
 & $B_{Bob,P_{separ}}$ & $O_{Bob,P_{separ}} $& $B_{Bob,P_{joint}} $& $O_{Bob,P_{joint}} $ \\  \hline
$B_{Ann,P_{separ}}$ &\begin{tabular}{c}   $(2;1)^{*,K=1}$  \\  $ \{\{1\},\{ 2\} \}$  \end{tabular} & \begin{tabular}{c} (0;0)  \\ $  \{\{1\},\{ 2\} \}$  \end{tabular} &\begin{tabular}{c}  (2;1) \\  $  \{\{1\},\{ 2\} \}$   \end{tabular} &  \begin{tabular}{c} (0;0)\\  $  \{\{1\},\{ 2\} \}$   \end{tabular} \\ \hline
$O_{Ann,P_{separ}}$ &  \begin{tabular}{c} (0;0) \\$  \{\{1\},\{ 2\} \}$ \end{tabular} & \begin{tabular}{c}  $(1;2)^{*,K=1}$ \\ $  \{\{1\},\{ 2\} \}$ \end{tabular}& \begin{tabular}{c}  (0;0) \\ $ \{\{1\},\{ 2\} \}$ \end{tabular} &  \begin{tabular}{c} (1;2) \\   $  \{\{1\},\{ 2\} \}$ \end{tabular} \\ \hline \hline
$B_{Ann,P_{joint}}$ &  \begin{tabular}{c} (2;1) \\ $  \{\{1\},\{ 2\} \}$ \end{tabular} &  \begin{tabular}{c} (0;0) \\ $ \{\{1\},\{ 2\} \}$ \end{tabular} &   
\begin{tabular}{c} 
$(2+\epsilon ;1+\epsilon )^{*,K=2} $
  \\ $ \{1, 2 \}$ 
  \end{tabular} & \begin{tabular}{c}  $ (\epsilon; \epsilon)$  \\  $  \{1,2\}$ \end{tabular} \\ \hline
$O_{Ann,P_{joint}}$ & \begin{tabular}{c}  (0;0) \\ $ \{\{1\},\{ 2\} \}$ \end{tabular} & \begin{tabular}{c}   (1;2) \\ $ \{\{1\},\{ 2\} \}$ \end{tabular}& 
\begin{tabular}{c}  $(\epsilon;\epsilon) $ \\ $ \{1,2\}$ \end{tabular} & 
\begin{tabular}{c}  $(1+\epsilon;2+\epsilon)^{*,K=2}$ \\  $ \{1, 2\}$ \end{tabular}  \\ \hline  
 \end{tabular}
\end{center}
\label{default2}
\end{table}%
 }}
 
 Formally the rules are: 
\begin{multline*}
 \mathcal{R}(K=1) \colon S(K=1) \mapsto S(\{ \{1 \} ,\{2 \} \}), \\  \forall s \in S_i(K=1)=
 \{O_{Ann,P_{separ}}, B_{Ann,P_{separ} } \} \times \{O_{Bob,P_{separ}},B_{Bob,P_{separ} }\}
  \end{multline*}
 and
 
 $$
 \mathcal{R}(K=2) \colon S(K=2) \mapsto \begin{cases}
  S(\{1,2 \}),  \\ \mbox{ if } s \in  \{O_{Ann,P_{separ}}, B_{Ann,P_{separ} } \} \times \{O_{Bob,P_{separ}},B_{Bob,P_{separ} }\}  
  \\ S(K=2) \setminus S(\{1,2 \} \}),  \ \ \ \ 
 \mbox{otherwise}
 \end{cases}
 $$
 
 The whole 
 Table  \ref{default2} corresponds to the game with $K=2$ where the game for $K=1$ is a  component.
  If Ann and Bob play the game with $K=1$ with a single coalition structure $\{\{ Ann\} , \{ Bob\} \}$, then the payoffs for this game are  in the two-by-two top-left corner of Table  \ref{default2}.   
 If Ann and Bob are together, then each obtains an additional payoff $\epsilon$,  and the corresponding cells make a two-by-two  bottom-right corner.
 
 Every game with $K=1$ and $K=2$ has only one    mixed strategies equilibrium and only one equilibrium partition.  Mixed strategies equilibrium for $K=1$ is described in every textbook.
 A change in $K$ from $K=1$ to $K=2$ results in a changes of an equilibrium strategy profile  and an equilibrium partition.
 
For $K=2$ the game  also has a mixed strategies equilibrium like for $K=1$. The difference is that now players have   an allocation in a different coalition structure in comparison to the case of $K=1$. Mixed strategies equilibrium  for Ann is: $\sigma^*(B_{Ann,P_{joint}})=(1+\epsilon)/(3+2\epsilon)$, 
 $\sigma^*(O_{Ann,P_{joint}})=(2+\epsilon)/(3+2\epsilon)$, $i \in \{\mbox{Ann}, \mbox{Bob} \}$, what results in an equilibrium partition $P_{joint}$, what is  the grand coalition. And there is  an equilibrium intra-coalition activity.
 
 The  presented game allows to construct intra-coalition externalities  from mixed strategies within one partition. Mixed intra-coalition externalities means that players are exposed to equilibrium fluctuations from strategic actions of another player.
 
A game as above  can not be constructed within any cooperative game equilibrium concept. It is impossible to construct Shapley value  (195) here even if there is one coalition: players have equilibrium mixed strategies inside it.

 \section{Stochastic games}
 Shapley  (1953) defined  stochastic games as " the play proceeds by steps from position to position, according to transition probabilities controlled jointly by the two player". 
This section demonstrates how this type of games with coalition structures as states of a game may appear in the constructed game.   The example  differs from example above as  a set of the equilibrium mixed strategies induces more than one equilibrium coalition structure.   We use a game, similar to corporate dinner game, but with identical players.
\subsection{Corporate lunch game}

There is a set of four identical players $N=\{A,B,C,D\}$. 
An individual strategy is a coalition structure, or an allocation of all players across tables during   lunch.   
A coalition structure is an allocation of   players over no more than four tables, where possibly empty tables do not matter. 
A rule of coalition formation is  a unanimous agreement to form a coalition. If  player  chooses a coalition, but other members of the coalition did not choose him, the player eats alone.

A player has preferences over coalition structures:  she/he prefers to eat with someone and  other players eat individually. If one eats alone he is hurt by a possible formed coalition of others.
 Coalition structures,  and payoff profiles  for the highest cases payoffs are presented  in Table  \ref{lunch}:
\begin{table}[htp]
\caption{Office lunch game: strategies and payoff profiles. Full set of equilibrium mixed strategies are indicated  only for player $A$}
 \begin{center}
\begin{tabular}{|c|c|c|c|}
\hline num&
Coalition structure &  Payoff profile $U_A,U_B,U_C,U_D$ & \begin{tabular}{c} \small{Coalition values } \\ \small{as in }  \\  \small{ cooperative game theory} \end{tabular}  \\
\hline
$1^*$& $\{ A,B\}, \{ C\}, \{D\}$:  $ \sigma^*_A= \sigma^*_B=1/3$ & (10,10,3,3)  & $20_{A,B}$, $3_{C}$, $3_{D}$ \\
$2^*$& $\{ A,C\}, \{ B\}, \{D\}$:  $ \sigma^*_A=\sigma^*_C=1/3$ & (10,3,10,3)  & $20_{A,C}$, $3_{B}$, $3_{D}$ \\
$3^*$& $\{ A,D\}, \{ C\}, \{B\}$:  $ \sigma^*_A=\sigma^*_D=1/3$ & (10,3,3,10)  & $20_{A,C}$, $3_{C}$, $3_{B}$ \\
4 & $\{ A\}, \{B\}, \{ C, D\}$ & (3,3,10,10)  & $3_{A}$, $3_{B}$, $29_{C,D}$\\
5 & $\{ A\}, \{D\}, \{ C, B\}$ & (3,10,10,3)  & $3_{A}$, $3_{B}$, $29_{C,D}$\\
6 & $\{ A\}, \{C\}, \{ B, D\}$ & (3,10,3,10)  & $3_{A}$, $3_{B}$, $29_{C,D}$\\
7& $\{ A\} ,\{B\}, \{ C\}, \{D\}$ & (3,3,3,3) & $3_{A}$, $3_{B}$, $3_{C}$, $3_{D}$\\
8 & \mbox{all other with $K=3,4$} & (0,0,0,0) &  $ = 0$\\
\hline
\end{tabular}
\end{center}
\label{lunch}
\end{table}%
 
 Payoffs in   Table \ref{lunch} are organized in the way that formation of  a coalition by other players deteriorates payoffs for the rest.  Formation of two 2-player coalitions does not improve payoffs for every player.  Outcomes of the game are coalition structures or states of a stochastic game. An increase  of $K=2$ to $K=3,4$ does not change an equilibrium in mixed strategies, hence we can speak about robustness of an equilibrium for $K=2$ to an increase in K. 

 It is clear that the game does not have an equilibrium in pure strategies. This is a Bayesian game, with a probability distribution of equilibrium mixed strategies. Equilibrium mixed strategies are indicated  only  for player A  in the first three lines of the Table \ref{lunch}.
 
 
 \subsection{A formal definition of a stochastic game of coalition structure formation}
 
 Let $\Gamma(K)$ be a non-cooperative game as defined above.
 
 \begin{definition}
 A game $\Gamma(K)$ is a stochastic game if a set of equilibrium partitions is bigger than two, $\# \Large( \{P^*\}(K)\Large) \ge 2$, where a state  is an equilibrium partition $P^* \in \{P^* \}(K)$.
 \end{definition}
 
 Studying properties of stochastic games with non-cooperative coalition structure formation are left for future.
  
  In the same way we can  construct a family of stochastic games, what is left for future.

 \section{No self-enforcement of an equilibrium and non-cooperative criterion for stability}

 Aumann (1990) used "stag and hare"   game to demonstrate absence of a  self-enforcement property of Nash equilibrium.  We use an extended version of the same game  to demonstrate how by modifying the game we can reach a focal point of the game, unavailable within standard Nash equilibrium. Then we    construct a  non-cooperative coalition structure stability criteria.
  
 There are two hunters $i=1,2$.
 If  players can hunt only individually, then $K=1$, and  the only partition is $P_{separ}=\{ \{ 1\}, \{ 2\}\}$.   An individual strategy  set  of $i$ is
$S_i(K=1)=\Big(({P_{separ}, hare}),({P_{separ}, stag}) \Big) $. For example, a strategy $s_i=({P_{separ}, hare})$ is interpreted as player $i$ chooses to  hunt alone for a hare.

A set of a corresponding strategies for the game with $K=1$ is $S(K=1)=S_1(K=1) \times S_1(K=1) $.  

For a game with $K=2$ every hunter  can choose either to hunt alone, in a coalition structure $P_{separ}=\{ \{ 1\}, \{ 2\}\}$, or  together, $P_{joint} = \{ 1,2\}$. For every hunting partition   a player  chooses a target for hunting: a hare or a stag, as in the game for $K=1$.
 A set  of strategies of $i$ is
 $$S_i(K=2)=\Big(({P_{separ}, hare}),({P_{separ}, stag}), ({P_{joint}, hare}),({P_{joint}, stag}) \Big),
 $$ where a strategy consists of two  terms. 
 A set of strategies of the game is a direct ( Cartesian) product, $S(K=2)=S_1(K=2) \times S_2(K=2)$.

  We do not rewrite the rules for coalition structure formation, as they are the same as  in the BoS game above. The difference is  a renaming of  strategies. 
  
Every player    knows, which game is played, either with $K=1$ or with $K=2$.  A case with uncertainty in   parameter $K$ is not addressed here and left for the future.

We assume that there is a unanimous coalition formation rule, i.e. hunters can hunt together   only if both choose to be together. However even if they both choose to hunt together, they may have a disagreement for whom to hunt.  
Payoffs for the both games $\Gamma(K=1)$ and $\Gamma(K=2)$ are presented in Table \ref{harestag}.

If hunters play a game with $K=1$ then a maximum achievable payoff is $(8,8)$, when each hunts individually for a hare. An an equilibrium strategy profile  is
$\Big((P_{separ}, hare),(P_{separ}, hare) \Big)$.
 In the game $\Gamma(K=1)$  the players  can not reach the  efficient outcome $(100,100)$  of the game $\Gamma(K=2)$, when both hunt for a stag together.
  This focal point (in terminology of Schelling) can  be  reached only within the game $\Gamma(K=2)$. 
Thus  the focal point is feasible in the game $\Gamma(K=2)$, but not in the game $\Gamma(K=1)$.

Self-enforcing property of the equilibrium  is that  both players realize individual gain from a change of a game from $K=1$ to $K=2$ and neither can deviate.
But players can not reach  the  outcome  $(100;100)$ without a change in a game.

In literature  a self-enforcement property of an equilibrium   is not well-defined, but intuitively it  depends  on  what   players think about willings of many others to deviate from an equilibrium.     
    \begin{table}[htp]
\caption{Expanded stag and hare game}
\begin{center}
\begin{tabular}{|l|c|c|c|c|}
\hline 
&\small{\mbox{\negthickspace  $P_{separ}$, hare}} &\small{\mbox{\negthickspace  $P_{separ}$, stag}} & \small{\negthickspace\mbox{ $P_{joint}$,hare}} & \small{\negthickspace \mbox{ $P_{joint}$, stag}} \\
\hline
\small{\mbox{$P_{separ}$,hare}} &  (8;8); $\{ \{ 1\},\{ 2\}\}$ & (8;0); $\{ \{ 1\},\{ 2\}\}$  & (8;8);$\{ \{ 1\},\{ 2\}\}$ & (8;0); $\{ \{ 1\},\{ 2\}\}$ \\
\small{\mbox{$P_{separ}$, stag} }  & (0;8); $\{ \{ 1\},\{ 2\}\}$ & (0;0); $\{ \{ 1\},\{ 2\}\}$& (0;8); $\{ \{ 1\},\{ 2\}\}$& (0;0); $\{ \{ 1\},\{ 2\}\}$\\
\hline
\small{\mbox{$P_{joint}$, hare} } & (8;8); $\{ \{ 1\},\{ 2\}\}$ & (8;0); $\{ \{ 1\},\{ 2\}\}$  & (4;4); $\{ 1, 2\}$ & (8;0); $\{ 1, 2\}$\\
 \small{\mbox{$P_{joint}$, stag}}   & (0;8); $\{ \{ 1\},\{ 2\}\}$& (0;0) ;$\{ \{ 1\},\{ 2\}\}$ & (0;8);$\{ 1, 2\}$ & $(100;100)^*$; $\{ 1, 2\}$\\
 \hline
\end{tabular}
\end{center}
\label{harestag}
\end{table}%

 If there is an uncertainty, which game is played,  either $\Gamma(K=1)$ or $\Gamma(K=2)$, then players will randomize between two strategies:
 $P_{separ}, hare$ and $P_{joint}, stag$. In this case 
 the game becomes a stochastic game as described above.
 
 \subsection{Criterion of  coalition structure (a partition) stability}
There is a    nested family  of games $${\Gamma} = \{\Gamma(K=1), \ldots, \Gamma(K), \ldots \Gamma(K=N) \} \colon \Gamma(K=1) \subset  \Gamma(K) \subset \Gamma(K=N).$$

It is   characterized by list of  equilibrium strategy profiles, $$\Big( \sigma^*(1), \ldots, \sigma^*(K), \ldots, \sigma^*(K=N)\Big),$$ where  $ \sigma^*(K)=( \sigma^*_i(K))_{i \in N}$ and by a list of  equilibrium  partitions
 $$\Big(\{ P^*\}(1), \ldots, \{ P^*\}(K),\ldots, \{ P^*\}(K=N) \Big),$$ $\{ P^*\}(K)  \subset \mathcal{P}(K)$. Every game $\Gamma (K)$  from a family $\Gamma $ has an equilibrium may be in mixed strategies.

The family of games has an equilibrium expected  payoff profiles:

$$
\Big( (EU^{\Gamma(1)}_i)^*_{i \in N}, \ldots,( EU^{\Gamma(K)}_i)^*_{i \in N}, \ldots,( EU^{\Gamma(K=N)}_i)^*_{i \in N} \Big) 
, $$ where $( EU^{\Gamma(K)}_i)^*_{i \in N} \equiv ( EU^{\Gamma(K)}_i (\sigma^*))_{i \in N}$. 

Let us take  a game $\Gamma(K_0) \in \Gamma$ with $\sigma^*(K_0)$ as an equilibrium mixed strategy set. The question is: what is a condition when   an equilibrium strategy profile  does not  change with an increase in a maximum  
coalition  size $K_0$?   
We can do this by comparing expected utility for  all players from   different games: $\Gamma(K_0), \ldots, \Gamma(K = N)$.


The  criterion is based on the fact that a set of mixed strategies should not change with an increase in a variety of available coalition structures.
 The criterium   is a sufficient criterium and defines a maximum coalition size, when an equilibrium   for smaller    coalition sizes still holds true.  
 \begin{definition}
\textbf{Partition stability criterion}   for a game $\Gamma(K_0)$ is a  maximum coalition size $K^*$ when an equilibrium  still holds true, i.e. for all $i \in N$ there is a maximum number $K^*$ such that
   \begin{enumerate}
   \item $$K^* =   \max_{K = K_0,\ldots,  N \atop{\Gamma(K_0) \ldots, \Gamma(K=N)}}  \Big\{ EU^{\Gamma(K_{0})}_i\Big(\sigma^*_i (K_{0}),
  \sigma^*_{-i}(K_{0})\Big) \ge \\  EU^{\Gamma(K)}_i\Big(\sigma^*_i (K),
  \sigma^*_{-i}(K)\Big)  \Big\},$$
  \item 
 $Dom \mbox{ }  \sigma^*(K^*)  = Dom \mbox{ }  \sigma^*(K_0)$
 \end{enumerate}
where $\sigma^*(K_0)$ is an equilibrium in the game $\Gamma(K_0)$,  $\sigma^*(K)$ is an equilibrium in a game $\Gamma(K)$, $K = K_0,\ldots,  N$, and  $Dom$  is a  domain of equilibrium mixed strategies set.

\end{definition} 

The definition is operational, it can be constructed directly from a definition of a family of games.
This definition guarantees stability of both payoffs and partition, and is a sufficient criterion of stability. Some applications may require weaker forms of the  criterion. 


Now it is clear that the statement of Aumann (1990) that Nash equilibrium is generally not self-enforcing is correct.  In the extended version of stag and hare game we have seen that an increase in $K$ changed an equilibrium. 
The same took place in a variation of Battle of Sexes game. 
However this did not happen in the Corporate  Dinner or the Corporate Lunch game.

The proposed criterion may serve as a measure of trust to an equilibrium or as a test for self-enforcing   of an equilibrium.
This criterion can be  applied to study opportunistic behavior  in coalition partitions. If  players  in a coalition $g$ in a game $\Gamma(K_1)$  have perfect cooperation, this does not mean that in a wider game  $\Gamma(K_2)$,  $K_1 < K_2$,  they will still cooperate. 

\section{Discussion}

The current paper demonstrates how  some standard non-cooperative games    can be constructed from the game offered by Levando (2016), i.e.  by including coalition structures into consideration. The proposed detail  for a  definition of a non-cooperative game  improves interpretability of    results and wides applications for non-cooperative games.
Using the example, the paper offers a way to construct  cooperation in  coalition formation on   self-interest fundamentals. The paper offers also a non-cooperative criterion  to measure stability of coalition structures using a definition of a non-cooperated games from a nested family.

The two accompanying papers demonstrate applications of the same model to network games (a generalization of Myerson, 1977) and to simple repeated games. 




\end{document}